\documentclass[preprint,3p,sort]{elsarticle}

\usepackage{amsfonts}
\usepackage{amsthm}
\usepackage{graphicx}
\usepackage{epstopdf}
\usepackage{algorithmic}
\usepackage{amsmath}
\usepackage{amssymb}
\usepackage{hyperref}
\usepackage[T1]{fontenc}
\usepackage{soul,color}
\usepackage{cases}
\newtheorem{theorem}{Theorem}[section]

\theoremstyle{remark}
\newtheorem{remark}[theorem]{Remark}
\theoremstyle{example}

\numberwithin{equation}{section}
\begin{document}


\begin{frontmatter}

%

\title{Uncertain Data in Initial Boundary Value Problems: \\ Impact on Short and Long Time Predictions}

\author[sweden,southafrica]{Jan Nordstr\"{o}m}
\cortext[secondcorrespondingauthor]{Corresponding author}
\ead{jan.nordstrom@liu.se}
\address[sweden]{Department of Mathematics, Link\"{o}ping University, SE-581 83 Link\"{o}ping, Sweden}
\address[southafrica]{Department of Mathematics and Applied Mathematics, University of Johannesburg, P.O. Box 524, Auckland Park 2006, Johannesburg, South Africa}

\begin{abstract}
We investigate the influence of uncertain data on solutions to initial boundary value problems. Uncertainty in the forcing function, initial conditions and boundary conditions are considered and we quantify their relative influence for short and long time calculations. It is shown that dissipative boundary conditions leading to energy bounds play a crucial role. For short time calculations, uncertainty in the initial data dominate. As time grows, the influence of initial data vanish exponentially fast. For longer time calculations, the uncertainty in the forcing function and boundary data dominate, as they grow in time. Errors due to the forcing function grows faster (linearly in time) than the ones due to the boundary data (grows as the square root of time). Roughly speaking, the results indicate that for short time calculations, the initial conditions are the most important, but for longer time calculations, focus should be on modelling efforts and boundary conditions.  Our findings have impact on predictions where similar mathematical and numerical techniques are used for both short and long times as for example 
in regional weather and climate predictions.

\end{abstract}



\end{frontmatter}


\section{Introduction}

Approximative solutions which constitute the output to initial boundary value problems (IBVPs) are typically generated by various forms of numerical schemes. The investigation of stability, accuracy and convergence on these approximations dominate in the numerical analysis litterature (see \cite{kreiss1970,kreiss1989initial,Gustafsson1978,gustafsson1995time,oliger1978,nordstrom2020,nordstrom_roadmap,nordstrom2005} for examples regarding linear problems and  \cite{godunov1961interesting,volpert1967,kruzkov1970,lax1973,harten1983,tadmor1984,Tadmor1987,Tadmor2003,nordstrom2021linear} for nonlinear ones). However, no output can be better than the quality of the input, which for IBVPs consists of forcing functions, initial data and boundary data.
This fact has historically attracted less interest and often the quality of the data is assumed to be high enough. However, that the quality of the input data often is problematic have been  realised in application areas such as aerospace \cite{doi:10.2514/1.3961,YAO2011450,YONDO201823}, nuclear physics \cite{Furnstahl_2015,Beane_2015}, oil prospecting \cite{CHRISTIE2006143,DOSTERT20083445,MONDAL2010241,10.2118/119139-PA}, weather forecasts \cite{doi:10.3402/tellusa.v65i0.21740,https://doi.org/10.1002/qj.3545,npg-10-211-2003} and climate predictions \cite{EstimatingtheUncertaintyinaRegionalClimateModelRelatedtoInitialandLateralBoundaryConditions,ATutorialonLateralBoundaryConditionsasaBasicandPotentiallySeriousLimitationtoRegionalNumericalWeatherPrediction,ImpactofInitialConditionsversusExternalForcinginDecadalClimatePredictionsASensitivityExperiment,SensitivityofTyphoonTrackPredictionsinaRegionalPredictionSystemtoInitialandLateralBoundaryConditions,https://doi.org/10.1029/2008JD010969} to name a few.  

One motivation for this work was the different views on external input found in applications areas where otherwise similar computational technology is used. In aeronautical flow investigations using computational fluid dynamics (CFD) for example, solutions (e.g. steady state and flutter solutions) are obtained by long time integration and the initial data is often ignored (being typically an arbitrary constant). The specifics of the solution is assumed to be given by the governing equations combined with its boundary conditions and the governing IBVP is essentially viewed as a boundary value problem (BVP).  This view is quite different from the one in numerical weather prediction (NWP), where the focus is on creating appropriate initial data (using so called data assimilation) and boundary conditions are seen as less important. The IBVP in this case is typically seen as an initial value problem (IVP).  Both in CFD and NWP various modelling efforts include uncertain parameter values. In CFD those are e.g. present in turbulence models  \cite{Wilcox} and shock treatments \cite{Engquist108369}. For NWP they are employed in models for cloud formation, rainfall, and various transport phenomena \cite{OCONNELL19963}. Appropriate choices of these uncertain parameters in the governing equations are necessary for successful predictions and can be described by uncertainty in the forcing function.

In many of the application areas above, so called uncertainty quantification (UQ) has been used to address the problem \cite{Ghanem200663,FERSON2004355,Crestaux20091161,SMITH13,Gardiner09,MR3328389}. 
The forward propagation procedure can roughly be described as: given a certain input with a probability distribution, how can its influence on the statistically distributed output-solution be estimated? The inverse procedure (data assimilation) can roughly be described as: how can the appropriate statistical input be arranged such that a reasonable correct statistical output-solution is obtained? The UQ approach is pragmatic, informative and rather technical in nature but sometimes disregard fundamental properties of the governing IBVP (although exceptions exist \cite{Nordstrom20151,Wahlsten2018192,Nordstrom2019463}).

In this paper, we will proceed differently and focus on the fundamental IBVP properties. We will employ previous knowledge about IBVPs used in so called error bounded schemes \cite{nordstrom2007error,kopriva2017error,nordstrom2018long}. These schemes use well posed boundary conditions that generate a damping term in the energy rate that subsequently lead to error bounds. We will in addition use the recent development for IBVPs in \cite{nordstrom2022linear-nonlinear,Nordstrom2022_Skew_Euler,Nordstrom2023_Nonlinear_BC_modif} where it is shown that 
an energy bound of both the nonlinear and linearised IBVP can be obtained if a skew-symmetric form of the governing equations 
is available. We will combine these techniques mentioned above and study the effects of errors in the forcing function, initial data and boundary data on the solution. 


The results from the continuous analysis roughly described above holds also for numerical calculations with stable schemes and sufficiently fine meshes. Such stable schemes can systematically be constructed by discretising the equations in space using summation-by-parts (SBP) operators \cite{svard2014review,fernandez2014review} which discretely mimic the integration-by-parts (IBP) procedure.  The dissipative boundary conditions that produce the damping term mentioned above can be inserted weakly using penalty terms as described  in \cite{Nordstrom2023_Nonlinear_BC_modif,nordstrom2017roadmap} which lead to a strongly stable schemes. The mimicking properties of schemes based on SBP operators lead directly to discrete error bounds similar to those found in the continuous formulation. We will not discuss the numerical approximation procedure in this paper, since that part is covered in numerous previous publications (by us and others), see the references above.
 
The remaining part of paper is organised as follows: In Section~\ref{sec:theory} we reiterate and complement the main findings in  \cite{nordstrom2022linear-nonlinear,Nordstrom2022_Skew_Euler,Nordstrom2023_Nonlinear_BC_modif,nordstrom2007error,kopriva2017error,nordstrom2018long} and outline the general procedure for obtaining energy bounds in terms of external data. The error analysis  is given in Section~\ref{error_estimates}. 
A summary is provided in Section~\ref{sec:conclusion}.

\section{The governing initial boundary value problem}\label{sec:theory}

We start by recapitulating the results in \cite{nordstrom2022linear-nonlinear,Nordstrom2022_Skew_Euler,Nordstrom2023_Nonlinear_BC_modif}. Consider the following general hyperbolic IBVP
\begin{equation}\label{eq:nonlin}
P U_t + (A_i(\bar U) U)_{x_i}+A^T_i(\bar U)U_{x_i}=F(\vec x,t),  \quad t \geq 0,  \quad  \vec x=(x_1,x_2,..,x_k) \in \Omega
\end{equation}
augmented with the initial condition $U(\vec x,0)=H(\vec x)$ in $\Omega$ and  the non-homogeneous boundary condition
\begin{equation}\label{eq:nonlin_BC}
L(\bar U) U = G(\vec x,t),  \quad t \geq 0,  \quad  \vec x=(x_1,x_2,..,x_k) \in  \partial\Omega.
\end{equation}
 In (\ref{eq:nonlin}) and in the rest of the paper, Einsteins summation convention with summation over repeated index  is used. The time-independent matrix  $P$ is symmetric positive definite and defines a scalar product $(U,V)_P= \int_{\Omega} U^T P V d\Omega$ and an energy norm $\|U\|^2_P=(U,U)_P$. We further require that the eigenvalues of $P$ (and hence also of $P^{-1}$)  are of order one. In (\ref{eq:nonlin_BC}), $F$ is a forcing function, $L$ is the boundary operator and $G$ the boundary data.  $F$, $G$ and $H$ is the external input data. We assume that both $U$ and $\bar U$ are smooth. The $n \times n$ matrices $A_i$ are smooth functions of the $n$ component vector $\bar U$, but otherwise arbitrary. Note that (\ref{eq:nonlin}) and (\ref{eq:nonlin_BC}) encapsulates both linear ($\bar U \neq U$) and nonlinear  ($\bar U=U$) problems. 




In \cite{Nordstrom2022_Skew_Euler,Nordstrom2023_Nonlinear_BC_modif} it was proved that the skew-symmetric form of (\ref{eq:nonlin}) with boundary conditions such that
\begin{equation}\label{1Dprimalstab_strong_extra}
\oint\limits_{\partial\Omega}U^T  (n_i A_i)   \\\ U \\\ ds = \oint\limits_{\partial\Omega} \frac{1}{2} U^T ((n_i A_i)  +(n_i A_i )^T) U \\\ ds \geq - \oint\limits_{\partial\Omega} G^TG \\\ ds
\end{equation}
lead to energy conservation and an energy bound.


\subsection{Modelling the effect of disturbed data}\label{3_energy_estimates} 
Consider the nonlinear problem (\ref{eq:nonlin}) with solution $U$ and disturbed data $F+\delta F$, $G+\delta G$ and $H+\delta H$. By subtracting the  non-disturbed problem with solution $V$ from the disturbed one and linearising using the technique in \cite{nordstrom2022linear-nonlinear} we find the evolution problem for the difference  $W=U-V$ to be
\begin{eqnarray}\label{eq:nonlin_diff}
P W_t + (A_i(U) W)_{x_i}+A^T_i(U)W_{x_i}&=&\delta  F(\vec x,t),  \quad t > 0,  \quad  \vec x=(x_1,x_2,..,x_k) \in \Omega \nonumber \\
L(U) W &=& \delta  G(\vec x,t),  \quad t >0,  \quad  \vec x=(x_1,x_2,..,x_k) \in  \partial\Omega  \\
W&=& \delta  H(\vec x), \quad \,\,\  t = 0,  \quad  \vec x=(x_1,x_2,..,x_k) \in  \Omega. \nonumber
\end{eqnarray}
By the development in \cite{Nordstrom2022_Skew_Euler,Nordstrom2023_Nonlinear_BC_modif} we know that $U$ in (\ref{eq:nonlin_diff}) is bounded by data with an appropriate boundary operator $L(U)$. 
We will investigate the influence of the disturbed data  $\delta F$, $\delta G$ and $\delta H$ on the deviation $W$ in  (\ref{eq:nonlin_diff}) by assuming that only one of the error sources are non-zero at a time, which enable us to assess their relative influence. The three types of errors in  (\ref{eq:nonlin_diff}) are: $\delta F$ which represent errors in modelling parameters, $\delta G$ which represent errors in boundary data, and finally $\delta H$ which represent errors in initial data.

\subsection{The dissipative boundary conditions}\label{Interpretation} 
Next, we recapitulate the crucial role of the boundary operator originally studied in  \cite{nordstrom2007error,kopriva2017error,nordstrom2018long}. 
The energy method applied to (\ref{eq:nonlin_diff}) yields
\begin{equation}\label{eq:boundaryPart1}
\frac{1}{2} \frac{d}{dt}\|W\|^2_P + \oint\limits_{\partial\Omega}W^T  (n_i A_i)  \\\ W \\\ ds= \int\limits_{\Omega}W^T  \delta  F d \Omega = (W, \delta  F)_I,
\end{equation}
where $(n_1,..,n_k)^T$ is the outward pointing unit normal and $I$ is the identity matrix. In (\ref{eq:boundaryPart1}), only the symmetric part of  $n_i A_i$ remains. Next we rotate the boundary term into the form
\begin{equation}\label{1Dprimalstab_trans}
\oint\limits_{\partial\Omega}W^T  (n_i A_i) W ds = \oint\limits_{\partial\Omega}C^T   \Lambda   \\\ C\\\ ds = \oint\limits_{\partial\Omega}(C^+)^T   \Lambda^+   \\\ C^+ + (C^-)^T   \Lambda^-   \\\ C^-\\\ ds,
\end{equation}
where $T^T  (n_i A_i) T =  \Lambda = diag( \lambda_i)$, $C = T^{-1} W$. In (\ref{1Dprimalstab_trans}), $\Lambda^+$ and  $\Lambda^-$ denote the positive and negative parts of $\Lambda$ respectively, while $C^+$  and $C^-$ denote the corresponding rotated variables.  We will use a dissipative boundary condition \cite{Nordstrom2023_Nonlinear_BC_modif} of the form  
\begin{equation}\label{Char_Bcond}
  \sqrt{|\Lambda^-|}C^-=\delta G.
\end{equation}

\begin{remark}\label{Sylvester}
For linear problems, the number of boundary conditions is equal to the number of eigenvalues of $ (n_i A_i)^S $ with the wrong sign  \cite {nordstrom2020}. Sylvester's Criterion \cite{horn2012}, show that the number of boundary conditions is equal to the number of $\lambda_i$  with the wrong sign if the rotation matrix  $T$ is non-singular.  In a nonlinear case it is more complicated since multiple forms of the boundary terms may exist since $\Lambda = \Lambda(C)$ \cite{nordstrom2022linear-nonlinear,Nordstrom2022_Skew_Euler,Nordstrom2023_Nonlinear_BC_modif,nordstrom2021linear}.  \end{remark}

\section{Error estimates due to the uncertainty in external data}\label{error_estimates}
In the upcoming estimates we will apply three relations. The first and second one are
\begin{equation}\label{time_norm_relation}
(W, \delta  F)_I  \leq \|W\|_P   \| \delta  F \|_{P^{-1}} \quad  \mbox{and}  \quad \frac{1}{2} \frac{d}{dt}\|W\|^2_P=\|W\|_P \frac{d}{dt}\|W\|_P.
\end{equation}
Following  \cite{nordstrom2007error,kopriva2017error,nordstrom2018long}, we next relate the outflow boundary terms to the $L_2$ norm of the solution as
\begin{equation}\label{Outflow_norm_relation}
\frac{ \oint\limits_{\partial\Omega}W^T  (n_i A_i) W ds}{\oint\limits_{\Omega}W^T  P W d \Omega} = \frac{ \oint\limits_{\partial\Omega} (C^+)^T   \Lambda^+  C^+ + (C^-)^T   \Lambda^-  C^-  ds}{ \int\limits_{\Omega}W^T  P W ds} \geq \frac{ \oint\limits_{\partial\Omega} (C^+)^T   \Lambda^+  C^+  ds}{\|W\|_P^2} = \eta(t).
\end{equation}
The relation $\eta(t)$ was shown in \cite{nordstrom2007error,nordstrom2018long} to lead to an integrating factor $\exp{(\theta(\xi,t))}$ in (\ref{eq:boundaryPart1}) where
\begin{equation}\label{int_factor}
 \theta(\xi,t)= \int_\xi^t \eta(\tau) d \tau
\geq \delta_0 (t-\xi) \quad  \mbox {with}  \quad \delta_0 > 0.
\end{equation}
The function  $\theta(\xi,t)$ is monotonically increasing in time since $\oint (C^+)^T   \Lambda^+  C^+  ds$, does not vanish for all time.

\subsection{Three different error estimates}\label{3_estimates}
Firstly we estimate the error or deviation for the case where $\delta F \neq 0, \delta G=0, \delta H=0$.  The relation (\ref{eq:boundaryPart1}) augmented with the homogeneous version of (\ref{Char_Bcond}), the relation (\ref{time_norm_relation}) and  (\ref{Outflow_norm_relation}) leads to 
\begin{equation}\label{DF_eq}
\frac{d}{dt}\|W\|_P + \eta(t)\|W\|_P   \leq \|\delta F\|_{P^{-1}}.
\end{equation}
The use of the integrating factor technique and the estimate (\ref{int_factor}) leads to 
\begin{equation}\label{DF_intermidiate}
\\|W\|_P   \leq \int_0^t   e^{-\theta(\tau,t)}  \|\delta F\|_{P^{-1}} d \tau  \leq \int_0^t   e^{-\delta_0 (t-\tau)} d \tau  (\|\delta F\|_{P^{-1}})_{\max\limits_{(0,t)}} \leq \frac{1-e^{-\delta_0 t }}{\delta_0} (\|\delta F\|_{P^{-1}})_{\max\limits_{(0,t)}}.
\end{equation}

Secondly we estimate the error or deviation for the case where $\delta F = 0, \delta G \neq 0, \delta H=0$.  The relation (\ref{1Dprimalstab_trans}) augmented with the non-homogeneous version of (\ref{Char_Bcond})
and the relation (\ref{Outflow_norm_relation})leads to 
\begin{equation}\label{DG_eq}
\frac{d}{dt}\|W\|_P^2 + 2 \eta(t)\|W\|_P^2   =  - 2\oint\limits_{\partial\Omega} (C^-)^T   \Lambda^-  C^-  ds =  2\oint\limits_{\partial\Omega}  (\delta G)^T \delta G ds=2 \|\delta G\|_{\partial\Omega}^2 .
\end{equation}
The use of the integrating factor technique and estimate (\ref{int_factor}) leads in this case to 
\begin{equation}\label{DG_intermidiate}
\\|W\|_P^2   \leq 2\int_0^t   e^{-2\theta(\tau,t)}  \|\delta G\|_{\partial\Omega}^2  d \tau  \leq 2 \int_0^t   e^{-2\delta_0 (t-\tau)} d \tau  (\|\delta G\|_{\partial\Omega}^2) _{\max\limits_{(0,t)}} \leq \frac{1-e^{-2\delta_0 t }}{\delta_0} (\|\delta G\|_{\partial\Omega}^2) _{\max\limits_{(0,t)}}.
\end{equation}

Thirdly we estimate the error or deviation for the case where $\delta F = 0, \delta G = 0, \delta H \neq 0$.  The relation (\ref{1Dprimalstab_trans}) augmented with the homogeneous version of (\ref{Char_Bcond}) and the relation (\ref{Outflow_norm_relation})leads to 
\begin{equation}\label{DH_eq}
\frac{d}{dt}\|W\|_P^2 + 2 \eta(t)\|W\|_P^2   = 0 .
\end{equation}
The use of the integrating factor technique and estimate (\ref{int_factor}) leads in this case directly to 
\begin{equation}\label{DH_intermidiate}
\\|W\|_P^2   \leq e^{-2\theta(0,t)}  \|\delta H\|_P^2   \leq e^{-2\delta_0 t}  \|\delta H\|_P^2 .
\end{equation}

\subsection{Effects on short and long time calculations}\label{Short_Long_Times}
For long times,  the estimates (\ref{DF_intermidiate}),  (\ref{DG_intermidiate}),  (\ref{DH_intermidiate}) implies that the errors in initial data decay exponentially. The errors stemming from
modelling and boundary data continue to grow and the ones from modelling grows faster than the errors due to boundary data.  For short times,  the estimates (\ref{DF_intermidiate}),  (\ref{DG_intermidiate}),  (\ref{DH_intermidiate})  lead respectively to the leading order approximations
\begin{equation}\label{Short_times}
\\|W\|_P   \propto t   (\|\delta F\|_{P^{-1}})_{\max\limits_{(0,t)}}, \quad  \|W\|_P   \propto \sqrt{t}  (\|\delta G\|_{\partial\Omega}^2) _{\max\limits_{(0,t)}},  \quad  \|W\|_P  \propto \|\delta H\|_P.
\end{equation}
The estimates (\ref{Short_times}) show that errors in initial data dominate for small times. 



\section{Summary}\label{sec:conclusion}
We have investigated the influence of uncertain data on solutions to initial boundary value problems. Uncertainty in the forcing function, initial data and boundary data have been considered and their relative influence for short and long time calculations have been assessed. 
For short time calculations, uncertainty in the initial data dominate. As time grows, the influence of initial data vanish exponentially fast. For long time calculations, the uncertainty in the forcing function and boundary data dominate, as they grow in time. Errors due to the forcing function grows faster (linearly in time) than the ones due to the boundary data (which grows as the square root of time). The results indicate that for short time calculations, the initial conditions are important, but for long time calculations, focus should be on modelling efforts and boundary conditions. Our results have impact on calculations where similar mathematical and numerical techniques are used for both short and long times, such as in regional weather and climate predictions.



\section*{Acknowledgments}

JN was supported by Vetenskapsr{\aa}det, Sweden [2021-05484 VR] and the University of Johannesburg.


\bibliographystyle{elsarticle-num}
\bibliography{References_Jan,References_andrew,References_Fredrik}

\end{document}